\documentclass{amsart}

\usepackage{amssymb}

\newtheorem{theorem}{Theorem}

\newtheorem{cor}{Corollary}

\theoremstyle{remark}
\newtheorem*{rem}{Remark}

\newtheorem*{question}{Question}

\numberwithin{equation}{section}
\newcommand{\Z}{\mathbb Z}
\newcommand{\C}{\mathbb C}
\newcommand{\R}{\mathbb R}
\newcommand{\N}{\mathbb N}

\newcommand{\restrict}{{\,\upharpoonright\,}}

\newcommand{\abs}[1]{\lvert#1\rvert}
\newcommand{\dabs}[1]{\lVert#1\rVert}

\def\itlog{\operatorname{itlog}}
\def\sing{\operatorname{sing}}
\def\re{\operatorname{Re}}
\def \bi{{\underline{i}}}
\def \bj{{\underline{j}}}
\def \br{{\underline{r}}}
\def \supp {\operatorname{supp}}
\def \wt {\operatorname{wt}}

\DeclareMathAlphabet{\mathbf}{OML}{cmm}{b}{it}

\begin{document}
\title[Julia's Equation]{Julia's Equation and Differential Transcendence}

\author{Matthias Aschenbrenner}
\address{Department of Mathematics, University of California, Los Angeles, CA 90095-1555, U.S.A.}
\email{matthias@math.ucla.edu}
\thanks{The first-named author was partially supported by National Science Foundation Grant DMS-0556197. 
Part of the work on this paper was done while he was visiting Humboldt-Universit\"at Berlin during the fall of 2010, and he thanks this institution and Andreas Baudisch for their hospitality. Both authors would like to thank the referee for an extensive list of suggestions which greatly helped to improve the paper.}

\author{Walter Bergweiler}
\address{Mathematisches Seminar, Christian-Albrechts-Universit\"at zu Kiel, Ludewig-Meyn-Str. 4, D-24098 Kiel, Germany}
\email{bergweiler@math.uni-kiel.de}

\begin{abstract}
We show that the iterative logarithm of each non-linear entire function is differentially transcendental over the ring of
entire functions, and we give a sufficient criterion for  such an iterative logarithm   to be differentially transcendental over the ring of convergent power series.
Our results apply, in particular, to the exponential generating function of a sequence arising from work of Shadrin and
Zvonkine on Hurwitz numbers.
\end{abstract}


\maketitle

\section{Introduction and Main Results}

\noindent
In 1871, Schr\"oder~\cite{Schroeder1871}
suggested to study the iteration of a meromorphic function $f$ by
using the functional equation
\begin{equation}\label{1a}
\phi(\lambda z)=f(\phi(z))
\end{equation}
that now bears his name.  If $f$ satisfies this equation,
then the compositional iterates~$f^n$ of $f$ satisfy $\phi(\lambda^n z)=f^n(\phi(z))$,
so in principle we have an ``explicit'' expression for the iterates
of $f$ in terms of $\phi$ and its inverse function.
Schr\"oder gave various examples, e.g.,  $\phi(z)=\tanh z$, $\lambda=2$ and 
$f(z)=2z/(1+z^2)$, as well as Jacobian elliptic functions~$\phi$
which satisfy~\eqref{1a} for certain rational functions $f$.

K\oe nigs~\cite{Koenigs1884} considered the case that $f$
is holomorphic in a neighborhood of a fixed point $\xi$ and showed
that if the {\it multiplier}
$\lambda=f'(\xi)$ of $f$ at $\xi$ satisfies $\lambda\neq 0$ and $|\lambda|\neq 1$,
then~\eqref{1a} has a unique solution $\phi$ holomorphic in a neighborhood
of $0$ such that $\phi(0)=\xi$ and $\phi'(0)=1$.
Poincar\'e~\cite[p.~318]{Poincare1890}
observed that if $|\lambda|>1$ and if $f$ is rational,
then $\phi$ extends to a function meromorphic in the plane,
and if $|\lambda|>1$ and $f$ is entire, then $\phi$ is entire.
Therefore the solution $\phi$ of~\eqref{1a} for $|\lambda|>1$  is also called the
\emph{Poincar\'e function} of $f$ at $\xi$.

Schr\"oder~\cite[p.~303]{Schroeder1871}
expressed the opinion that the functions $f$ whose iterates
can be determined using~\eqref{1a}
are of rather special type.
One may argue
that the results of K\oe nigs and Poincar\'e say the opposite,
but support for Schr\"oder's view is given by a result of Ritt~\cite{Ritt1926}
which implies that only very few Poincar\'e functions 
are elementary functions.
In order to state Ritt's result, we recall that a holomorphic
function or, more generally, a formal power series $g$, is said to be \emph{differentially algebraic} if it satisfies 
an algebraic differential equation; that is, an equation 
of the form 
$$P\big(z,g(z),g'(z),\dots,g^{(n)}(z)\big)=0$$ where $P$ is
a non-zero polynomial in $2+n$ indeterminates (for some $n$) with constant coefficients;
if $g$ is not differentially algebraic, then $g$ is called \emph{differentially transcendental}.
Ritt's result says that a polynomial with a 
differentially algebraic  Poincar\'e function is conjugate to a monomial,
a Chebychev polynomial, or the negative of a Chebychev polynomial,
the corresponding
Poincar\'e functions being the exponential or trigonometric functions.
For rational functions there are additional cases
arising from the multiplication theorems of elliptic functions.
Poincar\'e functions of transcendental entire functions are always
differentially transcendental~\cite{Bergweiler1995}.

A family of meromorphic functions is called \emph{coherent} (or \emph{uniformly differentially algebraic}\/) if there exists
an algebraic differential equation which is satisfied by all functions
in the family, and \emph{incoherent} otherwise.
Boshernitzan and Rubel~\cite[Theorem~6.1]{Boshernitzan1986}
showed that a Poincar\'e function of a rational or entire function $f$
is differentially algebraic if and only if the family of iterates of~$f$ is
coherent. Thus the above results about differential transcendence
of Poincar\'e functions can be rephrased as results about
incoherence of iterates.

We now turn to the case where the multiplier 
$\lambda=f'(\xi)$ of $f$ at its fixed point $\xi$
does not satisfy the conditions  $\lambda\neq 0$ and $|\lambda|\neq 1$
required for K\oe nigs' theorem.
If $|\lambda|=1$, but $\lambda$ is not a root of unity, Schr\"oder's
equation still has a formal power series solution. The question whether
this series converges is rather delicate and forms the subject matter
of famous results of Siegel, Brjuno and Yoccoz; 
see \cite[Section~11]{Milnor2006} for a discussion.
However, regardless of whether the series converges or not,
it is differentially  transcendental whenever~$f$ 
is a non-linear rational or entire function~\cite{Becker1995}.

If $\lambda=0$, then instead of Schr\"oder's equation one 
considers B\"ottcher's equation. Again the solutions are 
differentially transcendental except in special cases~\cite{Becker1995}.

Suppose now that $\lambda$ is a root of unity. 
In this case the fixed point $\xi$ is also called \emph{parabolic}.
Passing to an iterate
of $f$ we may assume that $\lambda=1$. Assuming without loss of 
generality that $\xi=0$ we then write $f$ in the form
\begin{equation}\label{1b}
f(z)=z+\sum_{k=p}^\infty f_kz^k \quad \text{($p\geq 2$, $f_k\in\C$ for $k\geq 2$, $f_p\neq 0$).}
\end{equation}
A basic result of complex dynamics, called the 
Leau-Fatou flower theorem~\cite{Fatou,Leau}, says that there are $p-1$ domains
$L_1,\dots,L_{p-1}$, 
called \emph{petals} of $f$, such that $f(L_j)\subseteq L_j$ and the restriction
$f^n\restrict L_j\to 0\in \partial L_j$ as $n\to\infty$,
for  $j=1,\dots,p-1$. (See \cite[Section~10]{Milnor2006}.)
Moreover, 
the Abel functional equation $\phi(z+1)=f(\phi(z))$ has 
a holomorpic solution $\phi_j$ mapping the right half-plane to~$L_j$.
The functions $\phi_j$ are again differentially transcendental~\cite{Becker1995}.

A way to describe the iteration of $f$ not only in the petals
but in a full neighborhood of $0$ is based on the functional equation
\begin{equation}\label{1c}
\phi(f(z))=f'(z) \phi(z)
\end{equation}
which is named after Julia (e.g., in~\cite[Sections 3.5B and 8.5A]{Kuczma1990})
  or Jabotinsky (e.g., in~\cite{Aczel1988}).
It has a unique formal power series solution
\begin{equation}\label{1c1}
\phi(z)=f_p z^p +\sum_{k=p+1}^\infty \phi_k z^k \quad \text{($\phi_k\in\C$ for $k\geq p+1$),}
\end{equation}
which is called the \emph{iterative logarithm} of $f$ and denoted here by $\itlog(f)$.
The name iterative logarithm, introduced by
\'Ecalle~(see \cite[p.~8]{Ecalle1974} or \cite{Ecalle1975}), is explained by
the identity 
\[
\itlog(f^n)=n\itlog(f)
\qquad\text{valid for all $n\in\N$.}\]
The general solution of~\eqref{1c} is given by $\phi=\alpha\itlog(f)$
where $\alpha\in\C$.

The series in~\eqref{1c1} converges only in exceptional cases.
For example, a result of  Erd{\H{o}}s  and Jabotinsky~\cite{Erdoes1960}
in combination with results of Baker~\cite{Baker1964}
and Szekeres~\cite{Szekeres1964} shows 
that the only functions $f$ meromorphic in $\C$ and of the form~\eqref{1b} 
for which the series in~\eqref{1c1} converges
in some neighborhood of $0$
are the functions $f(z)=z/(1-cz)$ 
where $c\in\C$, with $\itlog(f)(z)=cz^2$.  (However,  \'Ecalle \cite{Ecalle1976} has shown that the iterative logarithm of a function $f$ holomorphic in a neighborhood of~$0$ satisfying~\eqref{1b}  is always Borel summable.)

It follows from the results in \cite{Bergweiler1995,Boshernitzan1986,Ritt1926}
that the iterative logarithm $\itlog(f)$ of a non-linear rational or entire
function $f$ is differentially transcendental; cf.\ the remarks at the end of
section~\ref{rpp}.
This can be viewed as an indication that the coefficient sequence $(\phi_k)_{k>p}$ is very irregular:
If a formal power series $y=\sum_{k} y_kz^k\in\C[[z]]$ is differentially algebraic, then the coefficient sequence~$(y_k)$
satisfies a certain (in general, non-linear) kind of recurrence relation~\cite[pp.~186--194]{Mahler}. Of particular importance in combinatorial enumeration is the class of {\it $D$-finite}\/ (also called {\it holonomic}\/) power series \cite[Chapter~6]{Stanley}. These are the formal power series  whose coefficient sequence satisfies a homogeneous {\it linear}\/ recurrence relation of finite degree with polynomial coefficients; equivalently  \cite[Proposition~6.4.3]{Stanley} those which satisfy a non-trivial {\it linear}\/ differential equation over~$\C[z]$.

In this paper, we  show that for entire functions we have an even stronger irregularity result.
To formulate this result we need some terminology: Given a subring $R$ of the ring $\C[[z]]$ of formal power series over $\C$ which is closed under differentiation, we say that $\phi\in \C[[z]]$ is \emph{differentially transcendental over $R$} if $\phi$ does not satisfy a non-trivial polynomial equation in $\phi$ and its derivatives
with coefficients from~$R$. (Thus ``differentially transcendental'' is synonymous with
``differentially transcendental over~$\C[z]$.'')  

\begin{theorem} \label{thm1}
Let $f$ be a  non-linear entire function of the
form~\eqref{1b}. Then $\itlog(f)$ is differentially transcendental over
the ring of entire functions.
\end{theorem}

Under an additional hypothesis we can even show that $\itlog(f)$
is differentially transcendental over
the ring $\C\{z\}$ of power series with positive radius of
convergence.
In order to state this hypothesis,  for an entire 
function $f$ we denote by $\sing(f^{-1})$ 
the set of singularities of the inverse function of $f$; 
see~\cite[Section~4.3]{Bergweiler1993} for a discussion of their role in
complex dynamics.
The set $\sing(f^{-1})$ coincides with the set of critical and (finite) asymptotic 
values of~$f$. 
Here a point $w\in\C$ is called a \emph{critical value} if there exists $\xi\in\C$ such that $f'(\xi)=0$
and $f(\xi)=w$ while $w$ is called an \emph{asymptotic value} if there exists a curve 
$\gamma\colon [0,1)\to\C$ such that $\gamma(t)\to\infty$ and $f(\gamma(t))\to w$ 
as $t\to 1$.
If $f$ is a polynomial, then we only have to consider critical values, since
polynomials have no finite asymptotic values.

The \emph{Speiser class}~${\mathcal S}$ consists of all non-linear entire
functions $f$ for which  $\sing(f^{-1})$ is finite. 
It plays an important role in complex dynamics; 
cf.~\cite{Bergweiler1993,Eremenko1992}.

The maximal domains $U_j$ ($j=1,\dots,p-1$) containing the petals $L_j$  such that
$f(U_j)\subseteq U_j$ and $f^n\restrict U_j\to 0$ as $n\to\infty$ 
are called \emph{Leau domains} of~$f$. If $z\in U_j$, then
$f^n(z)\in L_j$ for large~$n$.
A classical result of Fatou (see \cite[Theorem~7]{Bergweiler1993} or
\cite[Theorem 10.15]{Milnor2006}) says that
$U_j\cap\sing(f^{-1})\neq\emptyset$ for all~$j=1,\dots,p-1$.

\begin{theorem} \label{thm2}
Let $f\in {\mathcal S}$ be 
of the form~\eqref{1b}. Denote by $U_1,\dots,U_{p-1}$ the associated Leau domains 
and suppose that 
\begin{equation}\label{1e}
\sing(f^{-1})\subseteq\{0\}\cup\bigcup_{j=1}^{p-1}U_j.
\end{equation}
Then $\itlog(f)$ is differentially transcendental over $\C\{z\}$.
\end{theorem}

Examples to which Theorem~\ref{thm2} applies are 
$f_1(z)=z+z^2$ and $f_2(z)=e^z-1$.
The function $f_1$ has only one critical point at $-1/2$ and
$f(-1/2)=-1/4$ is the corresponding critical value.
The function $f_2$ has the only asymptotic value~$-1$ and 
no critical values. Thus $\sing(f_1^{-1})=\{-1/4\}$ and
$\sing(f_2^{-1})=\{-1\}$. It follows from the result of Fatou mentioned
above, or by direct computation, that $f_1$ and $f_2$ satisfy the
hypothesis of Theorem~\ref{thm2}.

Other examples are $f_3(z)=\sin z$ with $\sing(f_3^{-1})=\{1,-1\}$
and two Leau domains at~$0$, one containing $1$ and one containing~$-1$,
and $f_4(z)=ze^z$ with $\sing(f_4^{-1})=\{0,-1/e\}$.

The results of \cite{Bergweiler2002,Buff2002} imply that
if  $\re a>3/4$, then both critical points of 
$f(z)=z+z^2+az^3$  are in the Leau domain at $0$.
Thus $f$ satisfies
the hypothesis of  Theorem~\ref{thm2} if $\re a>3/4$.
In fact, this even holds \cite[p.~ 277]{Buff2002} if
$\re a\geq 3/4-1/(2\log 3)$.

\medskip

Theorem~\ref{thm2} suggests the following open question.

\begin{question}
Let $f$ be any transcendental entire function of the form~\eqref{1b}. Is $\itlog(f)$ differentially transcendental over $\C\{z\}$?
\end{question}

The iterative logarithm 
$$\itlog(e^z-1)=\frac{1}{2}z^2-\frac{1}{12}z^3+\frac{1}{48}z^4-\frac{1}{180}z^5+\frac{11}{8640}z^6-\frac{1}{6720}z^7+\cdots$$ 
of $f(z)=e^z-1$ 
is of particular interest since it is the exponential generating function~(egf) of a sequence
$$0,0, 1, -\frac{1}{2}, \frac{1}{2}, -\frac{2}{3}, \frac{11}{12}, -\frac{3}{4}, -\frac{11}{6}, \frac{29}{4}, \frac{493}{12}, -\frac{2711}{6}, -\frac{12406}{15}, \frac{2636317}{60}, \dots$$
of rational numbers which recently arose in a conjecture made by Shadrin and Zvonkine~\cite{SZ} (and proved in~\cite{Aschenbrenner2012}) in connection with a generating series for Hurwitz numbers, and also in another context (ongoing joint work of the first-named author with van den Dries and van der Hoeven on asymptotic differential algebra~\cite{AvdDvdH2013}). 
By Theorem~\ref{thm2}, its egf $\itlog(e^z-1)$ is differentially transcendental over $\C\{z\}$.
We do not know whether the {\it ordinary}\/ generating function (ogf) of this sequence is differentially transcendental over~$\C[z]$, let alone over $\C\{z\}$. (See \cite{Klazar} for some differential transcendence results over~$\C\{z\}$ for ogf's of sequences of combinatorial origin.)
We also do not know whether the coefficients $\phi_k$ of the power series $\phi=\itlog(e^z-1)\in \mathbb Q[[z]]$ are non-zero for all~$k \geq  3$. (A computation with MAPLE showed that $\phi_k\neq 0$ for $k=3,\dots,300$.)
Some general results about the coefficient sequence $(\phi_k)$ in the case where $\phi=\itlog(f)\in \C[[z]]\setminus\C\{z\}$ can be found in \cite{Jabotinsky}.

\medskip

The idea of the proof of  Theorem~\ref{thm1} is as follows.
Assuming that $\itlog(f)$ is differentially algebraic over $\C\{z\}$, we 
start with a differential equation which is ``minimal'' in a certain sense; cf.~ \S \ref{ade}.
We then use the functional equation of the iterative logarithm, i.e., equation~\eqref{1c}, to obtain a differential
equation with meromorphic coefficients which is satisfied by $f$ and all its iterates.
This implies that a Poincar\'e function $\psi$ associated to $f$ also satisfies such a differential equation.
Using a result of Steinmetz (Theorem~\ref{steinm} in \S\ref{ste}) we deduce that $\psi$ actually satisfies an
algebraic differential equation with constant coefficients.
This contradicts results about differentially algebraic Poincar\'e functions due to Ritt and
the second author; see \S\ref{rpp}.

The idea  of the proof of Theorem~\ref{thm2} is to assume that $\itlog(f)$ satisfies a differential equation with
coefficients analytic in a neighborhood of $0$ and then use inverse branches of $f$ to continue these
coefficients analytically to the whole plane. The conclusion then follows from Theorem~\ref{thm1}.
Some further remarks on the proof of  Theorem~\ref{thm2}  are made immediately after the proof.

\subsection*{Conventions and notations}
Throughout the paper, $i$, $j$, $m$, $n$, $p$ range over the set $\N=\{0,1,2,\dots\}$ of natural numbers.

\section{Preliminaries}\label{prelims}

\noindent
In this section we first introduce some basic terminology concerning differential polynomials used later.
We then recall more basic facts on repelling periodic points and Poincar\'e functions, in addition to the ones already appearing in the introduction. In the last part of this section we state a theorem of Steinmetz which is at the heart of the proof of Theorem~\ref{thm1}.

\subsection{Algebraic differential equations}\label{ade}
Let $R$ be a {\it differential ring,} that is, a commutative ring (with~$1$) equipped with a {\it derivation} of $R$, i.e., a map $f\mapsto f'\colon R\to R$ which is additive  and satisfies the Leibniz Rule:
$$(f+g)'=f'+g',\quad (f\cdot g)'=f\cdot g'+f'\cdot g\qquad\text{ for all $f,g\in R$.}$$
We let $f\mapsto f^{(n)}$ denote the $n$th compositional iterate of $f\mapsto f'$.  A subring $S$ of $R$ which is closed under $f\mapsto f'$ is called a {\it differential subring} of $R$, and in this case $R$ is called a {\it differential ring extension} of~$S$. 
For any 
$(r+1)$-tuple $\bi =(i_0,\dots,i_r)$ of natural
numbers and an element $y$ in a differential ring extension of $R$, put 
$$y^{\bi} := y^{i_0}(y')^{i_1}\cdots (y^{(r)})^{i_r}.$$
Let $Y$ be a differential indeterminate over~$R$. Then $R\{Y\}$ denotes the ring of
differential polynomials in $Y$ over $R$ (not to be confused with the ring $\C\{z\}$ of convergent power series with complex
coefficients in the indeterminate $z$).
As ring, $R\{Y\}$
is just the polynomial ring $R[Y,Y',Y'',\dots]$ in the distinct indeterminates
$Y^{(n)}$ over $R$, where as usual we write $Y=Y^{(0)}$, $Y'=Y^{(1)}$, $Y''=Y^{(2)}$, etc. 
We consider $R\{Y\}$ as the differential ring whose derivation
extends the derivation of $R$ and satisfies $(Y^{(n)})'=Y^{(n+1)}$
for every $n$. 
For $P(Y)\in R\{Y\}$ and $y$ an element of a differential ring extension of~$R$, we let~$P(y)$ be the element of that extension obtained
by substituting $y,y',\dots$ for $Y,Y',\dots$ in $P$,
respectively. 
We say that an element $y$ of a differential ring extension of $R$
is {\it differentially algebraic over~$R$}  if there is some $P\in R\{Y\}$, $P\neq 0$, such that $P(y)=0$, and if $y$ is not differentially algebraic over $R$, then~$y$ is said to be {\it differentially transcendental over~$R$.}

For $P\in R\{Y\}$, the smallest $r\in\N$ such that 
$P\in R[Y,Y',\dots,Y^{(r)}]$ is called the
{\it order} of the differential polynomial $P$. 
Let $P\in R\{Y\}$ have order $r$, and let $\bi=(i_0,\dots,i_r)$ range over $\N^{1+r}$.
We call $Y^{\bi}$ a   \emph{monomial}, and denote by $P_{\bi}\in R$ the coefficient of $Y^{\bi}$ in $P$. Thus $P$ can be uniquely written as 
\[ P=\sum_{\bi} P_{\bi}\,Y^{\bi}, \label{P} \]
where the {\it support}\/ of $P$, defined by
$$\supp P := \big\{\bi \colon P_{\bi}\neq 0\big\},$$
is finite.  We say that a monomial $Y^{\bi}$ \emph{occurs} in $P$
if $\bi\in\supp P$.
We set
$$\abs{\bi}:=i_0+\cdots+i_r,\qquad \dabs{\bi}:=i_1+2i_2+\cdots+ri_r.$$
For  $P\neq 0$ we call
$$\deg(P)=\max_{\bi\in\supp P} \abs{\bi},\qquad \wt(P)=\max_{\bi\in\supp P} \dabs{\bi}$$
the {\it degree} of $P$ respectively the {\it weight} of $P$. We say that $P\neq 0$ is {\it homogeneous} if $\abs{\bi}=\deg(P)$ for every $\bi\in\supp P$ and {\it isobaric} if $\dabs{\bi}=\wt(P)$ for every $\bi\in\supp P$.

For $r,s\in\N$ with $r\leq s$ we identify each $\bi=(i_0,\dots,i_r)\in\N^{1+r}$ with the tuple $(i_0,\dots,i_r,0,\dots,0)\in\N^{1+s}$ and
thus view $\N^{1+r}$ as a subset of $\N^{1+s}$.
We set $\N^*:=\bigcup_{r\in\N} \N^{1+r}$ and
 order $\N^*$ anti-lexicographically; that is, for $\bi=(i_0,i_1,\dots)$ and $\bj=(j_0,j_1,\dots)\in\N^*$ we set
$$\bi < \bj\quad:\Longleftrightarrow\quad \text{there is  some $k\in\N$ with
$i_k<j_k$ and $i_{l}=j_{l}$ for $l\geq k+1$, }$$ 
and we set $\bi \leq \bj :\Longleftrightarrow \text{$\bi<\bj$ or $\bi=\bj$}$. It is easy to verify that $\leq$ is a well-ordering  of~$\N^*$, that is, $\leq$ is a linear ordering of $\N^*$, and every non-empty subset of $\N^*$ has a smallest element with respect to~$\leq$.
For $P\neq 0$ we let the {\it rank $\br=\br(P)$}\/ of $P$ be the  largest element of $\supp P$ with respect to~$\leq$. Below $\bi$, $\bj$ range over $\N^*$. For $\bi=(i_0,i_1,\dots)$ and $\bj=(j_0,j_1,\dots)$ we put $\bi+\bj=(i_0+j_0,i_1+j_1,\dots)$.

We view the ring $\C[[z]]$ of formal power series over $\C$ as a differential ring in the usual way (with derivation $\frac{d}{dz}$), and we work with two differential subrings of $\C[[z]]$:
the differential subring $\C\{z\}$ of $\C[[z]]$ consisting of the convergent power series, and
the smaller differential subring  $\C\{z\}_\infty$ of $\C[[z]]$ consisting of the (Taylor series at~$0$ of)
entire functions.
Theorem~\ref{thm1} says that the iterative logarithm of a non-linear entire function 
is differentially transcendental over~$\C\{z\}_\infty$, while Theorem~\ref{thm2}
says that---under the hypotheses made on~$f$--- the iterative logarithm of $f$ is differentially transcendental over~$\C\{z\}$. 

A differential field is a differential ring whose underlying ring happens to be a field.
Sometimes we find it convenient to work in the differential field of meromorphic functions
(which may be naturally identified with the fraction field of $\C\{z\}_\infty$, equipped with the unique derivation extending that of~$\C\{z\}_\infty$).
We note that a function is differentially algebraic over the field of meromorphic functions
if and only if it is differentially algebraic over the ring of entire functions.

\subsection{Repelling periodic points and Poincar\'e functions}\label{rpp}
Let $f$ be a non-linear entire (or rational) function.
A point $\xi\in\C$ is called a {\it periodic}\/ point of~$f$ if there exists some $p\geq 1$
such that $f^p(\xi)=\xi$; the smallest such $p$ is called the {\it period} of $\xi$.  
One calls a periodic point $\xi$ of $f$ with period $p$ \emph{repelling} if the \emph{multiplier
$\lambda=(f^p)'(\xi)$}\/ of $f$ at $\xi$ satisfies $|\lambda|>1$. 

The \emph{Julia set $J(f)$} of $f$ is the of all points in the plane
(or Riemann sphere)  where the iterates of $f$ do not form a normal family.
A standard result of complex dynamics says that $J(f)$ is the
closure of the set of repelling periodic points of $f$.
For rational functions this was already proved by Fatou and Julia, by different
methods (see
\cite[Section~14]{Milnor2006} for an exposition of both proofs),
for transcendental entire functions it is due to Baker~\cite{Baker1968} (see
\cite{Bargmann1999,BertelootDuval2000,Schwick} for simpler proofs).  
The Julia set of $f$ is always non-empty (in fact, a perfect set).

As mentioned in the introduction, results of K\oe nigs and
Poincar\'e say that if $\xi$ is a repelling periodic point of $f$ with period $p$ and multiplier $\lambda$,
then there exists a function $\psi$ holomorphic in a neighborhood
of $0$ such that $\psi(0)=\xi$, $\psi'(0)=1$ and
$\psi(\lambda z)=f^p(\psi(z))$, called the Poincar\'e function of $f$ at $\xi$. 
If~$f$ is rational, then $\psi$ is meromorphic in the plane, and
if~$f$ is entire, then so is~$\psi$.
Moreover, $\psi$ is given by (cf.\ \cite[p.~670]{Ritt1926})
\begin{equation}\label{k}
\psi(z)=\lim_{n\to\infty}f^{n p}(\xi+\lambda^{-n}z).
\end{equation}
Differentiating~\eqref{k} we also obtain
\[\psi^{(m)}(z)=\lim_{n\to\infty}\lambda^{-mn}(f^{np})^{(m)}(\xi+\lambda^{-n}z)\quad\text{for each $m$,}\]
hence
\begin{equation}\label{l}
\psi^{\bi}(z) = \lim_{n\to\infty}\lambda^{-\dabs{\bi}n}(f^{np})^{\bi}(\xi+\lambda^{-n}z)\quad\text{for each $\bi$}
\end{equation}
and thus
\begin{equation}\label{n}
(\psi')^{\bi}(z)=
\lim_{n\to\infty}\lambda^{-(\abs{\bi}+\dabs{\bi})n}\big((f^{np})'\big)^{\bi}(\xi+\lambda^{-n}z)\quad\text{for each $\bi$.}
\end{equation}
As mentioned in the introduction,
Ritt~\cite{Ritt1926} determined all differentially algebraic Poincar\'e functions
of rational functions.
His result shows in particular that rational functions with differentially algebraic
Poincar\'e functions have no parabolic fixed points, 
so there is no
iterative logarithm associated to these functions.
Moreover, 
it was shown in~\cite{Bergweiler1995} that Poincar\'e functions
to transcendental entire functions are differentially transcendental.
Combining this with Ritt's result we obtain the following.
\begin{theorem} \label{ritt} 
Let $f$ be a non-linear rational or entire function with a parabolic fixed point.
Then the Poincar\'e functions associated to the repelling fixed points of $f$ 
are all differentially transcendental.
\end{theorem} 
Together with the results of Boshernitzan and Rubel~\cite{Boshernitzan1986}
quoted earlier this implies the  following 
result already mentioned in the introduction:

\begin{cor}
Let $f$ be a non-linear rational or entire function.
Then
the iterative logarithm of $f$ at each parabolic fixed point of $f$
is differentially transcendental.
\end{cor}
\begin{proof}
Suppose $0$ is a parabolic fixed point of $f$; it is enough to show that then $\itlog(f)$ is differentially transcendental. Assume otherwise. Then by \cite[Theorem~6.4]{Boshernitzan1986} (see also \cite[Corollary~6.3]{Aschenbrenner2012}) there is a nonzero differential polynomial $P\in\C[z]\{Y\}$
such that $P(f^n)=0$ for all $n$. Let $\zeta\in\C$ be a repelling periodic point of $f$, with period $p$.  Replacing $f$ by $f^p$ we may assume that $p=1$, so $f(\zeta)=\zeta$.
Let $g(z):=f(z+\zeta)-\zeta$; then $0$ is a repelling fixed point of $g$, and with $Q:=P(Y+\zeta)$
we have $Q(g^n)=0$ for each $n$. Let $\psi$ be the Poincar\'e function of $g$ at~$0$.
By \cite[Theorem~6.1]{Boshernitzan1986}, $\psi^{-1}$ is differentially algebraic, hence so is~$\psi$, contradicting Theorem~\ref{ritt}.
\end{proof}

\subsection{A result of Steinmetz}
\label{ste}
The following result is due to Steinmetz~\cite[Satz~1]{Steinmetz1980}. We denote by
$T(r,f)$ the Nevanlinna characteristic of a meromorphic function $f$, and
as usual in Nevanlinna theory, 
$S(r,f)$ denotes any term satisfying
$S(r,f)=o\big(T(r,f)\big)$ as $r\to\infty$
outside some exceptional set of finite measure. See~\cite{Hayman1964}
as a reference for Nevanlinna theory.

\begin{theorem} \label{steinm} 
Let $F_0,F_1,\dots,F_m$
and $h_0,h_1,\dots,h_m$ be not identically vanishing meromorphic
functions
and let $g$ be a nonconstant entire function
such that
\[
F_0(g)h_0+ F_1(g)h_1+ \cdots + F_m(g)h_m=0.
\]
Suppose that there exists a positive  $K\in\R$ such that
\[
\sum_{j=0}^mT(r,h_j)\leq K\,T(r,g)+S(r,g).
\]
Then there exist polynomials
$P_0,P_1,\dots,P_m$ with constant coefficients, not all zero,  
such that
\[
P_0(g)h_0+ P_1(g)h_1+ \cdots + P_m(g)h_m=0 .
\]
\end{theorem}

\section{Proof of Theorem~\ref{thm1}}\label{pf of thm1}

\noindent
Let $f$ be a non-linear entire function as in \eqref{1b}, with iterative logarithm $\phi=\itlog(f)$.
Differentiation of~\eqref{1c} yields
\begin{equation}\label{e0}
\phi'(f)\cdot f' =f''\cdot \phi+f'\cdot \phi'=A_{01}(f') \cdot\phi + A_{11}(f')\cdot \phi'
\end{equation}
with $A_{01}(X)=X'$ and $A_{11}(X)=X$. 
Differentiating this equation, multiplying by $f'$  and substituting \eqref{e0}, one obtains
\begin{align*}
\phi''(f)\cdot(f')^3	&=\big(f'''f'-(f'')^2\big)\cdot \phi+f''f'\cdot \phi'+(f')^2\cdot\phi'' \\
						&= A_{02}(f') \cdot\phi + A_{12}(f')\cdot \phi' + A_{22}(f')\cdot \phi''
\end{align*}
with $A_{02}(X)=X''  X-(X')^2$, $A_{12}(X)=X'  X$  and $A_{22}(X)=X^2$. 
Induction yields the existence of differential polynomials 
$A_{ij}\in\Z\{X\}$ ($i\leq j$) in a differential indeterminate~$X$, 
independent of $f$, such that
\[
\phi^{(j)}(f)\cdot(f')^{2j-1} =\sum_{i=0}^j A_{ij}(f') \cdot\phi^{(i)}.
\]
Each $A_{ij}$ is homogeneous and isobaric, and if non-zero, 
of degree~$j$  and weight~\mbox{$j-i$}. (See \cite[Section~6.5]{Aschenbrenner2012}, where $H_{ij}(X)=A_{ij}(X')$.)
Moreover, $A_{jj}=X^j$ and the monomial of
highest rank occurring in $A_{0j}$ is $X^{(j)}X^{j-1}$.
For $\bj\in\N^*$ this yields
\begin{equation}\label{e}
\phi^{\bj}(f)\cdot(f')^{2\dabs{\bj}-\abs{\bj}} = 
\sum_{\substack{\bi\leq \bj\\ \abs{\bi}=\abs{\bj}}} B_{\bi,\bj}(f') \cdot\phi^{\bi} 
\end{equation}
with differential polynomials $B_{\bi,\bj}\in\Z\{X\}$ ($\bi\leq\bj$, $\abs{\bi}=\abs{\bj}$), independent of $f$. 
For $\bj=(j_0,\dots,j_r)\in\N^{1+r}$, we have
\[
B_{\bj,\bj} = \prod_{k=0}^r \left(A_{kk}\right)^{j_k}=X^{\dabs{\bj}}
\]
and 
\[
B_{(\abs{\bj}),\bj} = \prod_{k=0}^r \left(A_{0k}\right)^{j_k},
\]
so $B_{(\abs{\bj}),\bj}$ is homogeneous of degree $\dabs{\bj}$ and isobaric of weight $\dabs{\bj}$, and
the monomial of highest rank occurring in
$B_{(\abs{\bj}),\bj}$ is $X^{\bj}  X^{\dabs{\bj}-\abs{\bj}}$.
Note that for each $n$,~\eqref{1c} also holds with $f$ replaced by the iterate $F=f^n$ of $f$, that is,
\begin{equation}\label{c1}
\phi(F(z))=F'(z)\phi(z)\qquad\text{where $F=f^n$,}
\end{equation}
and so~\eqref{e} also holds with $f$ replaced by $F$.

Towards a contradiction assume now that $\phi=\itlog(f)$ is differentially algebraic 
over~$\C\{z\}_\infty$, that is, $\phi$ satisfies an equation
\begin{equation}\label{c}
P(\phi)=\sum_{\bi} P_{\bi} \phi^{\bi}=0
\end{equation} 
where $P=\sum_{\bi} P_{\bi} Y^{\bi}$ is a non-zero differential polynomial 
with entire coefficients~$P_{\bi}=P_{\bi}(z)$.
We assume that $P$ is chosen so that its rank $\br=\br(P)$ is  minimal.
Note that~$\dabs{\br}>0$, since otherwise \eqref{c} would show that $\phi$ is algebraic over $\C\{z\}$
and hence in~$\C\{z\}$ (since $\C\{z\}$ is algebraically closed in $\C[[z]]$ by Puiseux's Theorem, see
\cite[III, \S{}4]{Ruiz1993}),
contrary to the results in~\cite{Baker1964,Erdoes1960,Szekeres1964}
already quoted in the introduction, which say 
that the only functions $f$ meromorphic in $\C$  for which $\itlog(f)\in\C\{z\}$ are 
those of the form $f(z)=z/(1-cz)$.
Allowing the coefficients $P_{\bi}$ to be meromorphic, we may also assume that $P_{\br}=1$.

It follows from~\eqref{c}  and~\eqref{e} that
\[
\begin{aligned}
0
&=\sum_{\bj\leq \br} P_{\bj}(f) \cdot \phi^{\bj}(f) \cdot (f')^{2\dabs{\br}-\abs{\br}}\\
&=\sum_{\bj\leq \br}  P_{\bj}(f)\cdot (f')^{2\dabs{\br}-2\dabs{\bj}-\abs{\br}+\abs{\bj}}
\sum_{\substack{\bi\leq \bj\\ \abs{\bi}=\abs{\bj}}} B_{\bi,\bj}(f')\cdot \phi^{\bi}
\end{aligned}
\]
so that
\begin{equation}\label{f}
\sum_{\bi} \left(\sum_{\substack{\bi\leq\bj\leq\br \\ \abs{\bi}=\abs{\bj}}} P_{\bj}(f)\cdot (f')^{2\dabs{\br}-2\dabs{\bj}-\abs{\br}+\abs{\bj}}\cdot B_{\bi,\bj}(f')
 \right) \phi^{\bi} =0.
\end{equation}
It also follows from~\eqref{c} that
\begin{equation}\label{g}
\sum_{\bi} P_{\bi}\cdot (f')^{\dabs{\br}} \cdot \phi^{\bi} =0.
\end{equation}
In the last two equations, the coefficient of $\phi^{\br}$ is $(f')^{\dabs{\br}}$. 
By the minimality of $\br$ the two equations are thus equal. 
(We note that the exponent of $f'$ might actually be negative for some
terms on the left hand side of~\eqref{f}, but this does not affect the argument,
since we may multiply both equations by a sufficiently high power of $f'$. Similar adjustments are  tacitly made in what follows.)

Equating coefficients  in~\eqref{f} and~\eqref{g}
we obtain, for all $\bi< \br$, a (possibly trivial) differential equation for~$f$ with
meromorphic coefficients. We shall only consider the case that $\bi=(|r|)$ and we shall see, that
then the resulting differential equation for~$f'$ is non-trivial.
So we compare the coefficients of 
$\phi^{\bi}$ in~\eqref{f} and~\eqref{g} for $\bi=(|r|)$ and, putting $a=P_{(\abs{\br})}$,
we obtain
\begin{equation}\label{h}
\sum_{\substack{(\abs{\br})\leq \bj\leq \br\\ \abs{\bj}=\abs{\br}}} 
P_{\bj}(f) \cdot (f')^{2\dabs{\br}-2\dabs{\bj}}\cdot B_{(\abs{\bj}),\bj}(f')
=a\cdot (f')^{\dabs{\br}}.
\end{equation}
As noted before, $X^{\dabs{\bj}-\abs{\bj}}X^{\bj}$ is the  monomial of highest rank occurring in
$B_{(\abs{\bj}),\bj}$,
and thus the monomial of highest rank  in $X^{2\dabs{\br}-2\dabs{\bj}}\cdot B_{(\abs{\bj}),\bj}$ is $X^{2\dabs{\br}-\dabs{\bj}-\abs{\bj}}X^{\bj}$.
Hence among the monomials occurring in 
the differential polynomials  
$X^{2\dabs{\br}-\dabs{\bj}}B_{(\abs{\bj}),\bj}$
on the left side of~\eqref{h}
the one of maximal rank given by 
 $X^{\dabs{\br}-\abs{\br}}X^{\br}$,
 and it is contributed only by the term corresponding to~$\bj=\br$.
 Since $P_\br=1\neq 0$, we conclude that the differential equation~\eqref{h} is non-trivial.

Also, since  $B_{(\abs{\bj}),\bj}$ is homogeneous of degree $\dabs{\bj}$ and isobaric of weight~$\dabs{\bj}$,
each $\bi\in\supp X^{2\dabs{\br}-2\dabs{\bj}}B_{(\abs{\bj}),\bj}$ satisfies $\abs{\bi}+\dabs{\bi}=2\dabs{\br}$.
Thus, 
incorporating the terms $X^{2\dabs{\br}-2\dabs{\bj}}$ into the monomials
occurring in $B_{(\abs{\bj}),\bj}$, equation~\eqref{h} takes the form
\begin{equation}\label{j}
\sum_{\substack{\bi\leq (\dabs{\br}-\abs{\br})+\br \\ \abs{\bi}+\dabs{\bi}=2\dabs{\br}}} b_{\bi}(f) \cdot  (f')^{\bi}=a\cdot (f')^{\dabs{\br}}
\end{equation}
with meromorphic functions $b_{\bi}$, and $b_{(\dabs{\br}-\abs{\br})+\br}=1$.
Let $$I=\big\{\bi:\bi\leq (\dabs{\br}-\abs{\br})+\br,\  \abs{\bi}+\dabs{\bi}=2\dabs{\br} \big\}.$$
By the remarks following \eqref{e},
 equation \eqref{j} also holds for $f$ replaced by $F=f^n$, for each~$n$, so
\begin{equation}\label{j iterated}
\sum_{\bi\in I} b_{\bi}\big(F(z)\big) \cdot  \big(F'(z)\big)^{\bi}=a\cdot \big(F'(z)\big)^{\dabs{\br}}
\qquad\text{where $F=f^n$.}
\end{equation}
As noted in Section~\ref{rpp}, $f$ has repelling periodic points.
(Actually it was shown in~\cite{Bergweiler1991}
that every iterate of $f$ apart possibly from $f$ itself has repelling fixed points.)
Replacing~$f$ by an iterate, we may in fact assume that $f$ has a repelling fixed point~$\xi$.
Moreover, we may assume that $\xi$ is not a pole of~$a$.
With $\lambda=f'(\xi)$ we define the Poincar\'e function $\psi$ by~\eqref{k}.
From \eqref{l} and~\eqref{n} recall that
\begin{align*}
(\psi')^k(z) &= \lim_{n\to\infty}\lambda^{-kn}\big((f^{n})'\big)^k(\xi+\lambda^{-n}z),\\
(\psi')^{\bi}(z) &=
\lim_{n\to\infty}\lambda^{-(\abs{\bi}+\dabs{\bi})n}\big((f^{n})'\big)^{\bi}(\xi+\lambda^{-n}z)\quad\text{for each $\bi$.}
\end{align*}
We substitute $\xi+\lambda^{-n} z$ for $z$ in~\eqref{j iterated},
multiply both sides of the equation by~$\lambda^{-2\dabs{\br} n}$, take the limit as
$n\to\infty$, and using $\dabs{\br}>0$, obtain
\[
\sum_{\bi\in I} b_{\bi}(\psi)\cdot (\psi')^{\bi} =0.
\]
It is a standard result of Nevanlinna theory \cite[p.~56]{Hayman1964} that 
\[
T\big(r,\psi^{(k)}\big)\leq T(r,\psi)+S(r,\psi)
\]
for each $k$.
This implies that 
\[
\sum_{\bi\in I} T\big(r,(\psi')^{\bi}\big)\leq K \, T(r,\psi)+S(r,\psi)
\]
for some constant $K\in\R$. 
Theorem~\ref{steinm} now implies that there exist polynomials~$Q_{\bi}$ (${\bi\in I}$) with constant coefficients, not all zero, such that
\[
\sum_{\bi\in I} Q_{\bi}(\psi) \cdot (\psi')^{\bi} =0.
\]
Thus $\psi$ satisfies an algebraic differential equation
with constants coefficients, contradicting Theorem~\ref{ritt}.  \qed

\section{Proof of Theorem~\ref{thm2}}\label{pf of thm2}

\noindent
Suppose that $\phi=\itlog(f)$ satisfies an equation of the form~\eqref{c}
whose coefficients~$P_{\bi}$ are in $\C\{z\}$.
Again we assume the rank $\br=\br(P)$ of $P$ to be minimal.
We choose $\rho>0$ such that all $P_{\bi}$ are holomorphic in 
$D_\rho=\{z\colon  \abs{z}<\rho\}$.
Allowing the coefficients~$P_{\bi}$ to be meromorphic in $D_\rho$
we may assume that~$P_{\br}=1$.
We want to show that the $P_{\bi}$ are actually meromorphic in~$\C$,
thereby obtaining a contradiction to Theorem~\ref{thm1}.

Let $L_1,\dots,L_{p-1}$ be petals of $f$ associated to the fixed
point $0$ as stipulated in the Leau-Fatou theorem.
These petals $L_j$  can be chosen arbitrarily small, and thus we
may assume that their closures are contained in $D_\rho$.
As $f'(0)=1$, there exists a branch $\psi$ of the inverse function
of $f$ defined in a neighborhood of $0$ such that $\psi(0)=0$
and $\psi'(0)=1$. The Leau-Fatou theorem may also be applied
to $\psi$. We denote by $L_1',\dots,L_{p-1}'$ the petals for~$\psi$.
(These petals are also called repelling petals for~$f$.)

By \eqref{1e} there exists $n$ such that
\begin{equation}\label{4a0}
f^n\big(\sing(f^{-1})\big)\subseteq \{0\}\cup  \bigcup_{j=0}^{p-1} L_j\subseteq D_\rho
\end{equation}
and thus
\begin{equation}\label{4a}
f^m\big(\sing(f^{-1})\big)
\subseteq \{0\}\cup  \bigcup_{j=0}^{p-1} L_j\subseteq D_\rho
\quad \text{for all  $m\geq n$.}
\end{equation}
Again we put $F=f^n$
and, proceeding as in the proof of Theorem~\ref{thm1}, we find 
that the equations~\eqref{f} and~\eqref{g} are equal.
The coefficients of~\eqref{g} are defined in $D_\rho$ while
the coefficients of~\eqref{f}, with $f$ replaced
by~$F$, are defined in the component of $F^{-1}(D_\rho)$
that contains~$0$. 
We denote this component by~$V$.
So the germs of the coefficients $P_{\bi}$ at $0$ can be continued 
meromorphically to both $D_\rho$ and~$V$.

Actually, by passing to slightly smaller domains $L_j$ and $D_\rho$ 
if necessary, we may assume that these germs can be continued
meromorphically to a region containing the closure $\overline{V}$ of $V$.
Moreover, we may assume that 
$f^k\big(\sing(f^{-1})\big)\cap \partial D_\rho=\emptyset$ for all~$k$,
which implies that $\partial V$ consists of analytic curves.

By the choice of $n$ we have $\sing(f^{-1})\subseteq V$ and in fact
\begin{equation}\label{4b}
f^m\big(\sing(f^{-1})\big)\subseteq V \quad \text{for all $m$.}
\end{equation}
Also, we may choose the petals $L_j$ and $L_j'$
so small that $\overline{L_j}\subseteq V$ and 
$\overline{L_j'}\subseteq V$ for~$j=1,\dots,p-1$.

As mentioned, we want to
show that the germs of the coefficients $P_{\bi}$ at $0$ can be 
continued to functions meromorphic in $\C$.
By the Monodromy Theorem, it suffices to show that the germs can be continued
meromorphically along any curve in $\C$ starting in~$0$.
We may restrict here to curves which
intersect~$\partial V$ only finitely often.
For example, this follows since it suffices to consider continuation along
polygonal paths and since~$\partial V$ consists of analytic curves. 

We now show that it suffices to consider continuation
along those curves $\gamma\colon[0,1]\to\C$ for which there exists $t_1\in (0,1)$
such that $\gamma\big([0,t_1]\big)\subseteq \overline{V}$ 
while $\gamma\big([t_1,1]\big)\subseteq \C\setminus V$.
In fact, suppose that continuation along such curves is possible and let  
$\sigma\colon[0,1]\to\C$ be a curve 
such that 
$$\sigma\big([0,s_1]\big)\subseteq \overline{V},\ 
\sigma\big((s_1,s_2)\big)\subseteq \C\setminus \overline{V},\ 
\sigma\big([s_2,s_3]\big)\subseteq \overline{V}\quad (0<s_1<s_2\leq s_3\leq 1).$$ 
Then there exists a curve
$\tau\colon [s_1,s_2]\to\partial V$ satisfying
$\tau(s_1)=\sigma(s_1)$ and $\tau(s_2)=\sigma(s_2)$ 
which is homotopic to $\sigma\restrict [s_1,s_2]$ in $\C\setminus V$;
that is, there exists a continuous
function $\Gamma\colon [s_1,s_2]\times [0,1]\to \C\setminus V$ such that
$\Gamma(s,0)=\sigma(s)$ and $\Gamma(s,1)=\tau(s)$ for all $s\in [s_1,s_2]$ and
$\Gamma(s_1,t)=\sigma(s_1)$ and $\Gamma(s_2,t)=\sigma(s_2)$ 
for all $t\in [0,1]$.
Let $\sigma^*\colon[0,1]\to\C$ be defined by 
$$\sigma^*(t)=\begin{cases} \sigma(t)&\text{for $t\in[0,1]\setminus[s_1,s_2]$,} \\ 
\tau(t)&\text{for $t\in[s_1,s_2]$.}\end{cases}$$
By our assumption, meromorphic continuation is possible along
both ${\sigma\restrict[0,s_2]}$ and ${\sigma^*\restrict[0,s_2]}$.
Moreover, $\Gamma$ yields a homotopy from 
${\sigma\restrict[0,s_2]}$ to ${\sigma^*\restrict[0,s_2]}$ with the property 
that meromorphic continuation is possible along all curves in the homotopy.
Thus, by the Monodromy Theorem, meromorphic continuation 
along $\sigma\restrict[0,s_2]$ and $\sigma^*\restrict[0,s_2]$
leads to the same result.
Since $\sigma\restrict[s_2,s_3]=\sigma^*\restrict[s_2,s_3]$, 
meromorphic continuation 
along $\sigma\restrict[0,s_3]$ and $\sigma^*\restrict[0,s_3]$
also leads to the same result.

Now $\sigma^*$ has the property that $\sigma^*\big([0,s_3]\big)\subset \overline{V}$.
Starting with a path $\sigma\colon [0,1]\to \C$ which intersects
$\partial V$  in finitely many points 
$\sigma(s_1),\dots,\sigma(s_n)$, iteration of the above procedure yields
a path $\gamma\colon[0,1]\to\C$ 
such that continuation along $\sigma$ and $\gamma$ leads
to the same result,  and $\gamma$ has the additional property
that unless $\gamma\big([0,1]\big)\subseteq\overline{V}$,
there exists $t_1\in (0,1)$
such that $\gamma\big([0,t_1]\big)\subseteq \overline{V}$ 
and $\gamma\big([t_1,1]\big)\subseteq \C\setminus V$.
It thus suffices to consider curves $\gamma\colon [0,1]\to\C$
with $\gamma(0)=0$ for which such a $t_1$ exists.
We may also assume that $\gamma\restrict [0,t_1]$ is injective.

Let now $\gamma$ be such a curve.
We have to show that the germs of the $P_{\bi}$ at $0$ can be continued 
meromorphically along~$\gamma$.
In order to do so, we may deform the part of $\gamma$ which is in $V$,
as long as it stays in~$V$.
Thus we may 
choose $\gamma$ such that 
$\gamma(t)\in L_1'$ for $t\in(0,t_0]$ with some $t_0\in(0,t_1)$,
but
\begin{equation}\label{4c}
\gamma\big([0,1]\big)\cap \bigcup_{j=0}^{p-1} L_j = \emptyset.
\end{equation}
Using~\eqref{4a}, \eqref{4b} and~\eqref{4c}
and noting that $\sing(f^{-1})$ is finite by hypothesis, we may in fact
assume that
\begin{equation}\label{4d}
\gamma\big((0,1]\big)
\cap 
\overline{\bigcup_{l=0}^\infty f^l\big(\sing(f^{-1})\big)}
= \emptyset.
\end{equation}
This implies that branches of the inverse functions of the iterates
of $f$ defined in a neighborhood of $0$ can be continued analytically
along~$\gamma$. In particular, for each $m$ we may continue along~$\gamma$
the branch $\psi_{m,0}$
of the inverse function of $f^m$,
defined in some neigborhood $U_0$ of~$0$,
which is given by $\psi_{m,0}(0)=0$.
Thus for each $t\in (0,1]$,
there exists a neighborhood $U_t$ of $\gamma(t)$,
a holomorphic function $\psi_{m,t}\colon U_t\to\C$
and $\delta>0$ such that whenever $|s-t|<\delta$,
then $\gamma(s)\in U_t$ and  $\psi_{m,s}(z)=\psi_{m,t}(z)$
for all $z$ in some neighborhood of $\gamma(s)$.
Moreover, while~$U_0$ depends on $m$,
it follows from~\eqref{4d} that the domains $U_t$ may be 
chosen independent of $m$ for $t\in (0,1]$. For example, we may choose $U_t$ as the
largest disk around $\gamma(t)$ which does not intersect
$\overline{\bigcup_{l=0}^\infty f^l\big(\sing(f^{-1})\big)}$.

It also follows from~\eqref{4d} that there exists 
a simply connected domain $\Omega$ containing $\gamma\big((0,t_0)\big)$
such that
\[
\Omega
\cap 
\overline{\bigcup_{l=0}^\infty f^l\big(\sing(f^{-1})\big)}=\emptyset.
\]
Thus all $\psi_{m,0}$ may be continued to functions holomorphic in~$\Omega$ 
and we may in fact assume that $\Omega\subseteq U_0$ for all~$m$.
By~\cite[Theorem~9.2.1]{Beardon1991} the $\psi_{m,0}$ form a normal family
in $\Omega$.
Since $\psi_{m,0}\restrict (L_1'\cap\Omega)\to 0$ as $m\to\infty$ 
we deduce that in fact $\psi_{m,0}\restrict{\Omega}\to 0$ as~$m\to\infty$.
This implies that $\psi_{m,t}\to 0$ as $m\to\infty$ for all $t\in (0,1]$,
locally uniformly in the domains $U_t$ where they are defined.
Altogether we see that if $m$ is sufficiently large, then
$\psi_{m,t}\big(\gamma(t)\big)\in D_\rho$ for all $t\in [0,1]$.
For the curve $\sigma\colon [0,1]\to \C$ defined by 
$\sigma(t)=\psi_{m,t}\big(\gamma(t)\big)$ we thus have
$\sigma\big([0,1]\big)\subseteq D_\rho$.

Next we note that~\eqref{c1} also holds for negative $n$,
with a negative exponent standing for the branch of the inverse function
of the appropriate iterate of $f$ 
which fixes~$0$.
Thus
\[
\phi\big(\psi_{m,0}(z)\big)=\psi_{m,0}'(z) \,\phi(z)
\]
for $z$ in a neighborhood of $0$.
Using this instead of~\eqref{c1} we obtain~\eqref{f} and~\eqref{g}
with~$f$ replaced by $\psi_{m,0}$.
As before we find that these equations are equal.
For the equation corresponding to~\eqref{g} the coefficients are
meromorphic in $D_\rho$.
In particular, 
since $\sigma\big([0,1]\big)\subseteq D_\rho$,
the germs of the $P_{\bi}$ at $0$ can be continued meromorphically along~$\sigma$.
Noting that $\sigma(t)=\psi_{m,t}\big(\gamma(t)\big)$,
we deduce that the germs of the functions~$P_{\bi}(\psi_{m,0})$ at~$0$
can be continued meromorphically along~$\gamma$.
Since the 
coefficients of the equation corresponding to~\eqref{f} are built from
the $P_{\bi}(\psi_{m,0})$ and  from differential polynomials in~$\psi_{m,0}$,
these coefficients can also be continued meromorphically along~$\gamma$.
As the equations corresponding to~\eqref{f} and~\eqref{g} are equal,
 we see that
the~$P_{\bi}$ can  be continued meromorphically along~$\gamma$. \qed
 
\medskip

The basic idea of the above proof appears in a paper of Lewin~\cite{Lewin1965}
who proved that $\itlog(f)\notin\C\{z\}$ for $f=e^z-1$. 
Assuming that $\itlog(f)$ is holomorphic in $D_\rho$ but has
a singularity $\zeta\in\partial D_\rho$, it is shown there by elementary
estimates that $w_1=f(\zeta)\in D_\rho$ or that there exists $w_2\in D_\rho$
with $f(w_2)=\zeta$. These points $w_j$ are also singularities of $\itlog(f)$,
leading to a contradiction. Note that $w_2=\psi(\zeta)$ for some branch
$\psi$ of the inverse of~$f$. 
The proof of Theorem~\ref{thm2} also uses the idea that given 
$\zeta\in\partial D_\rho$ there exists~$m$ such that
$f^m(\zeta)\in D_\rho$ or $\psi_m(\zeta)\in D_\rho$
for some branch $\psi_m$ of the inverse function of $f^m$.
However, in this more general setting
we have to be careful about the domain where this branch
of the inverse can be defined.

\begin{rem}
The Eremenko-Lyubich class $\mathcal B$ is defined as the class of all non-linear entire
functions  for which $\sing(f^{-1})$ is bounded. The proof of Theorem~\ref{thm2} shows
that instead of demanding that $f\in{\mathcal S}$ it suffices to assume that $f\in{\mathcal B}$
and that there exists~$n$ such that \eqref{4a0} holds. This is equivalent to saying 
that on $\sing(f^{-1})$ the iterates of~$f$ converge uniformly to~$0$.
An example of a function to which this remark applies is given by 
$f(z)=(\sin^2 z)/z$.
\end{rem}

\end{document}